\documentclass[11pt,twoside]{amsart}
\usepackage{float}
\usepackage[utf8]{inputenc}
\usepackage[english]{babel}
\usepackage[T1]{fontenc}
\usepackage{mathtools}
\usepackage{enumitem}
\usepackage{bbm}
\usepackage{cancel}
\usepackage{tabularx}
\usepackage{graphicx}
\usepackage{tikz}
\usepackage[all]{xy}
\usepackage{tcolorbox}
\usepackage{amsmath,
    amsfonts,
    amssymb,
    mathrsfs,
    stackrel}
\usepackage{comment}
\usepackage{bookmark}
\usepackage{hyperref}
\usepackage{cleveref}
\hypersetup{colorlinks,
    citecolor=blue,
    linktocpage}
\usepackage[normalem]{ulem}

\DeclareFontFamily{U}{wncy}{}
\DeclareFontShape{U}{wncy}{m}{n}{<->wncyr10}{}
\DeclareSymbolFont{mcy}{U}{wncy}{m}{n}
\DeclareMathSymbol{\Sha}{\mathord}{mcy}{"58}

\newcommand\cA{{\mathcal A}}

\newcommand\bF{{\mathbb F}}

\newcommand\bN{{\mathbb N}}

\newcommand\bZ{{\mathbb Z}}

\newcommand\F{{\mathbb F}}

\newcommand\wC{{\widetilde C}}

\newcommand{\set}[1]{\left\{ #1 \right\}}
\newcommand{\va}[1]{\left| #1 \right|}

\DeclareMathOperator{\Aut}{Aut}
\DeclareMathOperator{\Sol}{Sol}
\DeclareMathOperator{\Alb}{Alb}

\newtheorem{theorem}{Theorem}[section]

\newtheorem{corollary}[theorem]{Corollary}
\newtheorem{lemma}[theorem]{Lemma}

\newtheorem{proposition}[theorem]{Proposition}
\newtheorem{fact}[theorem]{Fact}

\theoremstyle{definition}

\newtheorem{definition}[theorem]{Definition}
\newtheorem{example}[theorem]{Example}

\newtheorem{notation}[theorem]{Notation}

\newtheorem{question}[theorem]{Question}
\newtheorem{remark}[theorem]{Remark}

\title[IHC for $\cA_6$]{Optimality of the 
Prym--Tyurin construction
for $\cA_6$}
\author[Engel]{Philip Engel}
\address{Department of Mathematics, Statistics, 
and Computer Science, University of Illinois 
Chicago (UIC),
851 S Morgan St, Chicago, IL 60607, USA}
\email{pengel@uic.edu}
\author[de Gaay Fortman]{Olivier de Gaay Fortman}
\address{Department of Mathematics, 
Utrecht University, 
Budapestlaan 6, 3584 CD Utrecht, The Netherlands.}
\email{a.o.d.degaayfortman@uu.nl}
\author[Schreieder]{Stefan Schreieder}
\address{Leibniz University Hannover, 
Institute of Algebraic Geometry, 
Welfengarten 1, 30167 Hannover, Germany.}
\email{schreieder@math.uni-hannover.de}

\date{\today}

\usepackage[maxnames=5, backend=bibtex, style=alphabetic]{biblatex}
\begin{document}

\subjclass[2020]{Primary 14K10, 05C25, 14C25; 
Secondary 14H40, 05B35, 14C30}
\keywords{Integral Hodge conjecture, abelian varieties, 
graph invariants, Prym--Tyurin varieties, 
algebraic cycles, regular matroids}

\begin{abstract} 
We prove that on a very general 
principally polarized abelian $6$-fold, 
the smallest multiple of the minimal curve
class which can be represented by an 
algebraic cycle is $6$.
\end{abstract}

\setcounter{tocdepth}{1}
\maketitle 

%\tableofcontents

\section{Introduction}  
This paper
is a sequel to 
\cite{companion}, where we gave counterexamples
to the integral Hodge conjecture for abelian varieties,
arising from certain combinatorial invariants of 
regular matroids.

It was classically known by Riemann that a general
principally polarized abelian variety 
$(X,\Theta)\in \cA_g$ of dimension $g=1,2,3$ 
is the Jacobian of smooth
genus $g$ curve $C$. Thus, 
the integral Hodge class $\Theta^{g-1}/(g-1)!$
is represented by an algebraic cycle;
the Abel--Jacobi image of $C$.
For $g=4,5$, the general 
member $(X,\Theta)\in \cA_g$ was shown by 
Wirtinger in 1895 to
be the Prym variety of an \'etale double 
cover $\wC\to C$ of a curve $C$ of genus $g+1$.
Thus $2\cdot \Theta^{g-1}/(g-1)!$
is represented by an algebraic curve; the image
of $\wC$ under its natural map to the Prym variety. 

In fact, the multiple $2$ is optimal, for very
general principally polarized abelian varieties
of dimension $g=4,5$. By \cite[Thm.~1.1]{companion},
the minimal algebraic multiple of 
the minimal curve class
$\Theta^{g-1}/(g-1)!$ on a very general 
$(X,\Theta)\in \cA_g$
is divisible by $2$, for any $g\geq 4$.

\begin{question}
Let $(X,\Theta)\in \cA_g$ be very general. 
What is the smallest positive integer $N_g\in \bN$ 
for which $N_g\cdot \Theta^{g-1}/(g-1)!$ 
is algebraic?
\end{question}

A simple observation is that $N_g \,{\mid}\, (g-1)!$
since the algebraic multiples of $\Theta^{g-1}/(g-1)!$
form an ideal in $\bZ$, and $\Theta^{g-1}$ is visibly
algebraic. Furthermore $N_g\,{\mid}\, N_{g+1}$ follows 
by specializing to a product of an elliptic curve 
with a very general principally polarized 
abelian variety of lower dimension,
see \cite[Sec.~8.3, Proof of Thm.~1.1]{companion}. 
The purpose of this paper is to add
one entry to the list of known values of $N_g$. 
In tabular form:

 \begin{table}[H]
 \begin{tabular}{|c|c|c|c|c|c|c|c|c|}
 \hline
 $g$ & $1$ & $2$ & $3$ & $4$ & $5$ & $6$ \\
 \hline
 $N_g$ & $1$ & $1$ & $1$ & $2$ & $2$ & $6$ \\
 \hline
 \end{tabular} 
 \end{table}

\begin{theorem}\label{main}
Let $(X,\Theta)\in \cA_6$ be a 
very general principally polarized
abelian $6$-fold. The smallest positive 
multiple $N_6$
of the minimal curve class 
$\Theta^5/5!$ which can be represented by an 
algebraic cycle is exactly $6$.
\end{theorem}

By work of 
Alexeev--Donagi--Farkas--Izadi--Ortega 
\cite{ADFIO}, it is known that $N_6\,{\mid}\, 6$. 
To prove Theorem \ref{main}, it thus suffices to 
show that also  $6\,{\mid}\,N_6$. 
For this, we consider a certain positive 
integer-valued combinatorial invariant $d(G)$ of a 
graph $G$, introduced in \cite{companion}. In the 
case $G=K_{g+1}$ of the complete graph on $g+1$ 
vertices, the results in {\it loc.cit.}~imply 
that this invariant satisfies
$d(K_{g+1})  \,{\mid}\, N_g$, 
see Theorem \ref{thm:divides} below. 
Moreover, by \cite[Theorem 1.1]{companion}, 
we know that  $2 \,{\mid}\, N_6$.
Hence to prove that $6 \,{\mid}\, N_6$ 
it suffices to show $3  \,{\mid}\, d(K_7)$. 
With the initial version
of the program we wrote for computing such 
invariants, checking that $3\,{\mid}\, d(K_7)$ 
was computationally infeasible. 
Here, we develop some theoretical
tools, which allow us to ``average'' with respect
to certain symmetries of the graph $K_7$, 
in turn reducing the computation to a manageable
size. See Remark \ref{computation-size} for 
further details.

Prym--Tyurin constructions
\cite{welters, kanev, birkenhake-lange, salomon}
explicitly produce families
of abelian varieties $A$ on which some
multiple of the minimal curve class,
called the {\it exponent},
is the class of a curve $C \subset A$.
Upper bounds for $N_g$ arise
from  Prym--Tyurin constructions that 
dominate $\cA_g$. It is unknown whether
for large $g$, the minimal exponent 
of such a Prym--Tyurin construction 
is less than $(g-1)!$.
Conversely,
our invariants of regular matroids
give lower bounds for $N_g$.
For $g \leq 6$, we are lucky enough
to achieve equality between
the upper and lower bounds. 

As discussed in 
\cite[Prop.~8.10]{companion}, 
we can (independently of this paper)
show that $N_7 \geq 6$ by
computing $6\,{\mid}\, d(K_{3,5})$. 
But arriving at an exact value
of $N_7$ will likely require
(1) techniques that give obstructions
to the algebraicity of prime {\it power} 
multiples of the minimal class, and (2) 
a new Prym--Tyurin construction 
of small exponent, dominating $\cA_7$.

\begin{question}
What is the value of $N_7$? 
\end{question}

\subsection{Acknowledgments}
PE was partially supported by NSF grant DMS-2441240.
OdGF and StS have received funding from the European Research Council (ERC) under the
European Union’s Horizon 2020 research and innovation programme under grant agreement
N\textsuperscript{\underline{o}} 948066 (ERC-StG RationAlgic). 
OdGF has also received funding from the ERC Consolidator Grant FourSurf N\textsuperscript{\underline{o}} 101087365. 
The research was partly conducted in the framework of the DFG-funded research training group RTG 2965: From Geometry to Numbers, Project number 512730679. The Gemini HPC cluster of Utrecht University was used to perform computer calculations.

\section{Matroids and their solutions}\label{gen}

\subsection{Solutions of matroids in graphs} 
\label{sec:solutions-albanese}
There are many equivalent formulations
of the notion of a regular matroid, such
as a matroid realizable over any field, or
a matroid realizable by a totally unimodular
matrix. For the purposes of this note, we declare:

\begin{definition}
A {\it regular matroid} $R$ on a 
finite {\it ground set} $S$ 
is a matroid realized by a map 
$S\to U^*={\rm Hom}(U,\bZ)$, $s \mapsto x_s$, 
where $U$ is a free abelian group of finite rank, 
for which:
\begin{enumerate}
\item the $x_s$ generate $U^*$, and
\item any subset of the $x_s$ generates a 
saturated sublattice of $U^*$.
\end{enumerate}
\end{definition}

We call $g = \mathrm{rank}(U)$ the {\it rank} of $R$. 
By condition (2), the base changed map
$S\to U^*\otimes \bF$ realizes
the same matroid over any field $\bF$.
The map 
$S\to U^*$ defines canonically a map
$\bZ^S \to U^*$
sending $ e_s \mapsto x_s$. Dually,
we have an injective 
map $U\hookrightarrow (\bZ^S)^*\simeq \bZ^S$. 

Let $H$ denote an oriented graph and let 
$E(H)$ denote the set of oriented edges. 
An {\it $S$-coloring} of $H$ 
is a partition $E(H)=\sqcup_{s\in S}E_s(H)$ of the
edges, indexed by the set $S$.
We say that an edge has color $s$ if $e\in E_s(H)$.
Following \cite[Def.~5.8]{companion}: 

\begin{definition}\label{def:sol} 
A {\it $\F_p$-solution of $R$ in 
an $S$-colored graph $H$}
%over the finite field $\bF_p$ 
is the data $b = (b_s)_{s\in S}$ of
a collection of $1$-chains $b_s\in C_1(H,\bF_p)$, 
for $s\in S$, for which 
\begin{enumerate}
\item[(A)] 
$b_s$ is supported on edges $E_s(H)\subset E(H)$ 
of color $s\in S$, and
\item[(B)] whenever
$\sum_{s\in S} a_se_s \in U \subset \bZ^S$,
$a_s\in \bZ$, we have that
$\sum_{s\in S} a_s b_s \in H_1(H,\bF_p)$ 
is a cycle.
\end{enumerate}\end{definition}

In \cite{companion}, we considered the more 
general notion of \emph{$\Lambda$-solution} of $R$, for a $\bZ_{(p)}$-algebra $\Lambda$; 
in this paper, all solutions 
are %over $\F_p$ 
$\bF_p$-solutions. We use the 
term ``solution'' to mean ``$\F_p$-solution'' 
when no confusion can arise.
Since $H$ is oriented, we have a canonical identification $C_1(H,\F_p)=\F_p^{E(H)}$.
Since $E(H)=\sqcup_{s\in S}E_s(H)$, the sum map 
$\F_p^{E_s(H)}\to \F_p$ which sends a vector to the 
sum of its coefficients thus induces a color profile map 
$\lambda\colon C_1(H,\F_p)\to \F_p^{S}$.

\begin{definition} The {\it color profile} 
of a solution of $R$ in $H$ is the vector 
$\lambda(b)\coloneq \lambda(\sum_{s\in S} b_s) \in \bF_p^{S}.$
We say that the color profile is \emph{constant} if 
$\lambda(b) = (c, \dotsc, c) \in \bF_p^S$ for some $c \in \bF_p$.
\end{definition}

If $R$ admits a $\bF_p$-solution in some graph $H$, 
then it in fact admits a solution with the same color profile
in a certain ``universal graph'' called the 
{\it mod $p$ Albanese graph} 
${\rm Alb}_p(R)$, see \cite[Sec.~5]{companion} 
(also notated ${\rm Alb}_{p,1}(R)$ in {\it loc.cit.}). 
We recall its construction:

\begin{definition}\label{def:alb}
The vertices of ${\rm Alb}_p(R)$ 
are the elements of the $\bF_p$-vector space 
$\bF_p^S/U_{\bF_p}.$ The
edges of ${\rm Alb}_p(R)$ are given by connecting 
$v\rightarrow v + \overline{e}_s$ for any vertex 
$v\in \bF_p^S/U_{\bF_p}$
and $\overline{e}_s$ the image
of any basis vector $e_s\in \bF_p^S$. 
In this way, ${\rm Alb}_p(R)$ is a directed graph, 
which is naturally $S$-colored, by 
declaring the edge $v\to v+ \overline{e}_s$ 
to have color $s\in S$.
\end{definition}

The {\it Cayley graph} of a group with a 
distinguished set of generators is the graph 
whose vertices are the elements
of the group and whose directed edges are given by (left) multiplication of the generators. The basis
vectors of $\bZ^S$ 
mapping to the quotient
$\bZ^S/U$ realize
the {\it dual matroid} $R^*$, defined
on the same ground set $S$. 
Thus, in concise terms: ${\rm Alb}_p(R)$ is the Cayley
graph of the mod $p$ realization 
$S \to \bF^S_p/U_{\bF_p}$ 
of the dual matroid: 
${\rm Alb}_p(R) = {\rm Cay}_p(R^*)$ 
with respect to
the defining vectors of the dual matroid. 

\begin{remark} All of the above
definitions make sense over $\bZ$, at the expense of 
the graph ${\rm Alb}_\bZ(R)={\rm Cay}_\bZ(R^*)$ 
being infinite. \end{remark}

\begin{notation}\label{sol-sub} 
We use the term {\it solution of $R$ mod $p$}, 
or simply, \emph{solution of $R$}, 
without reference to a graph, to mean an
$\bF_p$-solution of $R$ in ${\rm Alb}_p(R)$.
The solutions of $R$ mod $p$ form an $\bF_p$-vector space, which we denote by ${\rm Sol}_p(R)$. 
By Definition \ref{def:sol}(A), the supports
of the $1$-chains $(b_s)_{s\in S}$ are disjoint.
So we can view 
$\Sol_p(R) \subset C_1(\Alb_p(R), \bF_p)$ 
as a subspace
of the space of $1$-chains of its mod $p$ Albanese graph, 
by taking the sum $\sum_{s\in S} b_s$.
\end{notation}

\begin{remark}
The graph $\Alb_p(R)$ satisfies a 
certain universal property 
\cite[Prop.~5.7]{companion}:
Any connected, oriented $S$-colored graph $H$ 
whose color profile map 
$\lambda\colon C_1(H,\F_p)\to \F_p^S$ satisfies 
$\lambda(H_1(H,\F_p))\subset U_{\F_p}\subset \F_p^S$
admits, for any choice of vertex $v_0\in V(H)$, a 
canonical map of oriented $S$-colored graphs 
${\rm alb}\colon H\to \Alb_p(R)$
which does not contract any edge.
In particular, the pushforward of a solution of 
$R$ in $H$ yields a solution in $\Alb_p(R)$ 
with the same color profile.
\end{remark}

\subsection{Graphic matroids and the minimal class}

Let $G$ be an oriented graph with edge set $E(G)= E$. We will
specialize the discussion of the previous section
as follows:
\begin{enumerate}
\item $S=E(G)=E$ is the set of edges of $G$.
\item $U \to \bZ^S$ is the inclusion
$B^1(G,\bZ)\to C^1(G,\bZ) = \bZ^{E}$ of 
$1$-coboundaries into $1$-cochains.
\item $\bZ^S \to U^\ast$ is the quotient map $\bZ^E = C_1(G,\bZ) \to  C_1(G,\bZ)/H_1(G,\bZ) = B^1(G,\bZ)^\ast$.
\end{enumerate}

Here, we observe that $H_1(G,\bZ)$ embeds naturally into 
$C_1(G,\bZ)\simeq \bZ^E$ by viewing a cycle as a 
$\bZ$-linear combination of oriented
edges $i\in E(G)$. 

For an edge $i \in E(G)$, let $e_i \in \bZ^{E(G)}$ 
denote the associated basis element.

\begin{definition} The {\it graphic matroid} $M(G)$
is the matroid with ground set the edges
$E=E(G)$, realized by the map $i\mapsto x_i \coloneqq 
\overline{e}_i \in \bZ^E/H_1(G,\bZ)$. \end{definition}

To define the dual matroid of $M(G)$, we first dualize 
the map $C_1(G,\bZ) \to C_1(G,\bZ)/H_1(G,\bZ)$,
giving the map 
$B^1(G,\bZ) \to C^1(G,\bZ)$. 
Then the dual matroid of the graphic matroid, 
the {\it cographic matroid}
$M^*(G)$, is realized by the map
$$E(G) \to H^1(G,\bZ) = 
C^1(G,\bZ)/B^1(G,\bZ)$$
which assigns to the 
$i$th edge the linear functional 
$e_i^* \in H^1(G,\bZ) = H_1(G,\bZ)^*$ 
whose value
on any $1$-cycle is the coefficient of the
edge $e_i$ in it. 

Applying Definition \ref{def:alb}
and the ensuing discussion,
the Albanese graph ${\rm Alb}_p(M(G))$
of the graphic matroid $M(G)$ is the 
Cayley graph of $H^1(G,\bF_p)$ 
under the generators $e_i^*\in H^1(G,\bF_p)$
defining the
cographic matroid $M^*(G)$.

Consider a graphic matroid $M(G)$. Its rank 
is $|E(G)|-h^1(G) = |V(G)| - h^0(G).$
The number of connected components of $G$ can always be arranged
to equal $1$, by identifying a vertex on each connected
component to a single vertex. This does not change the edge
set or the first homology, and so leaves the 
graphic and cographic
matroids the same. We always arrange for
$G$ to have one connected component, 
so that the rank of $M(G)$ is $|V(G)|-1$. 

\begin{definition}
A regular matroid 
$R$ with ground set $S$ is 
{\it simple} if no $x_s=0 \in U^\ast$ and no 
$x_s=\pm x_t \in U^\ast$ for $s\neq t$.
\end{definition}

In terms of graphs, $M(G)$ is simple if $G$ has no loops
or parallel edges. Thus, when $M(G)$ is simple, 
$G$ is a minor of the
complete graph on $|V(G)|$ vertices 
$K_{|V(G)|}$ because $G$ arises by
deleting some subset of the edges of $K_{|V(G)|}$. 
We deduce
the following (very simple) fact:

\begin{fact}\label{fact1} 
The complete graph $K_{g+1}$ defines the unique
rank $g$ graphic matroid $M(K_{g+1})$
that is maximal, under the ``minors'' 
partial ordering of simple graphic 
matroids of rank $g$.
\end{fact}

\begin{definition} Let $d(G)$ be the product of all primes
$p\in \bN$ for which $M(G)$ admits {\bf no} mod $p$ solution 
whose color profile $(c,\dots,c)$ is nonzero and 
equal to a fixed constant $c\in \bF_p^*$. \end{definition}

It follows \cite[Rem.~8.5]{companion} that
$\textstyle  
d(G)\leq \prod_{p<g\textrm{ prime}} p$.
If a regular matroid $R$ has a mod $p$ 
solution with  constant, nonzero color profile,
then all minors of $R$ have a mod $p$ 
solution with the same property, see
\cite[Prop.~7.2]{companion}. 
Combined with Fact \ref{fact1}, we deduce 

  \begin{proposition} \label{div} For all $g\in \bN$, the following holds:
  \begin{enumerate}
      \item for any connected, simple graph $G$ on $g+1$  vertices, $d(G)\,{\mid}\, d(K_{g+1})$;
      \item 
 $d(K_g)\,{\mid}\, d(K_{g+1})$.
  \end{enumerate}
 \end{proposition}

\begin{theorem}[\cite{companion}] \label{thm:divides} 
The integer $d(K_{g+1})$ divides 
the minimal 
multiple $N_g$ of the minimal class 
$\Theta^{g-1}/(g-1)!$
which is algebraic on 
a very general  $(X,\Theta)\in \cA_g$.
\end{theorem}

\begin{proof} 
This follows from \cite[Thm.~8.10]{companion} 
and \cite[Lem.~2.8]{companion}. 
\end{proof}

We can view $d(K_{g+1})$ as the ``graphic contribution''
to the matroidal obstruction of \cite{companion}.
In terms of matroidal monodromy cones of
degenerations \cite[Def.~4.6]{survey}, the relevant
cone is ``Voronoi's principal domain of the first
type'' \cite[Ex.~6.7]{survey}, and the base
of the corresponding
degeneration of abelian varieties
dominates $\cA_g$.

\subsection{Solutions of graphic matroids}
We now analyze in more detail the linear 
algebraic set-up of $\bF_p$-solutions
for the graphic matroid $M(G)$, associated 
to an oriented graph $G$. 
As noted, ${\rm Alb}_p(M(G))={\rm Cay}_p(M^*(G))$ and thus,
the vertex set of ${\rm Alb}_p(M(G))$ is $H^1(G,\bF_p)$
while the $E(G)$-colored edges 
are of the form $v \to v+e_i^*$.
So we have an isomorphism 
\begin{align} \label{isom1}
C_1({\rm Alb}_p(M(G)),\bF_p) \simeq \bF_p[H^1(G,\bF_p)]\otimes C^1(G,\bF_p)
\end{align}
where $\bF_p[\cdot]$ denotes the group algebra.
This isomorphism sends the edge 
$v\to v+e_i^*$ to $[v]\otimes e_i^*$.
Then the boundary map 
$\partial\colon 
C_1({\rm Alb}_p(M(G)),\bF_p)\to 
C_0({\rm Alb}_p(M(G)),\bF_p)$
is identified with
%Here $\partial$ is defined to be 
\begin{align}\label{partial}
\begin{aligned}
\partial \colon \bF_p[H^1(G,\bF_p)]\otimes C^1(G,\bF_p) 
& \to \bF_p[H^1(G,\bF_p)], \\
[v] \otimes e_i^* 
& \mapsto [v+e_i^*]-[v].
\end{aligned}\end{align}
Here, we view $C^1(G,\bF_p) = (\bF_p^{E})^\ast$, 
via the orientation of $G$.

\begin{proposition}\label{solution-rephrase}
Let $G$ be an oriented
graph with edges $E=E(G)$. 
The space of $\bF_p$-solutions $b = (b_i)_{i\in E}
\in {\rm Sol}_p(M(G))$
of the graphic matroid $M(G)$ is 
identified
with space of the maps
$$C^1(G,\bF_p)\xrightarrow{\beta} 
\bF_p[H^1(G,\bF_p)]\otimes C^1(G,\bF_p),$$
satisfying the following two conditions:
\begin{enumerate}
\item[{\rm (A)}] $\beta(e_i^*) \in -\otimes e_i^*$
for all $i\in E$, and 
\item[{\rm (B)}] $\beta(B^1(G,\bF_p))\subset \ker(\partial)$,
\end{enumerate}
where $\partial$ is defined
in (\ref{partial}).
Furthermore, define a map
\begin{align*} \lambda\colon
\bF_p[H^1(G,\bF_p)]\otimes
C^1(G,\bF_p)&\rightarrow C^1(G,\bF_p), \\
[v]\otimes e_i^*&\mapsto e_i^*. \end{align*}
Then $\lambda \circ \beta$ is a diagonal matrix 
$\oplus_{i\in E}\, \bF_p e_i^\ast \to 
\oplus_{i\in E}\,\bF_p e_i^\ast$
and the color profile of a solution 
$b\leftrightarrow \beta$ is 
$\lambda(\beta(e_i^\ast))_{i \in E} \in C^1(G,\bF_p)$, 
via the isomorphism $C^1(G,\bF_p) \simeq \bF_p^E$. 
In particular, the condition
that the color profile be
constant is the same as saying
$\lambda\circ \beta = c\cdot {\rm Id}_{C^1(G,\bF_p)}$
for some $c\in \bF_p$.
\end{proposition}

\begin{proof}
Definition \ref{def:sol}(A) is easily
seen to translate
into condition (A) in the statement
of the proposition, by taking $\beta(e_i^*)\coloneqq 
b_i$ under the isomorphism (\ref{isom1}).
Then, Definition \ref{def:sol}(B) exactly
says that, in case of a graphic matroid, $\beta$ sends 
$U = B^1(G,\bF_p)\subset C^1(G,\bF_p)$
into $H_1({\rm Alb}_p(M(G)),\bF_p)=\ker(\partial)$.

Furthermore, $(\lambda\circ \beta)(e_i^*)$ computes
the $i$-th entry of the color profile because
it is the sum of the coefficients of 
$b_i = \beta(e_i^*)$ in $\bF_p$.\end{proof}

\subsection{Symmetrizing solutions} 
We now consider ``symmetrization'' of solutions.
Let $G$ be an oriented graph with edges $E=E(G)$ and 
let ${\rm Aut}(G)$ 
be its automorphism group, not necessarily 
preserving the chosen orientation of the edges of $G$. 
Then ${\rm Aut}(G)$ acts canonically
on the spaces $C^1(G,\bF_p)$ and 
$H^1(G,\bF_p)$ via pullback of simplicial
cochains and cohomology.
It also acts on the $E$-colored graph
${\rm Alb}_p(M(G))$, i.e.~the Cayley
graph of $H^1(G,\bF_p)$ with respect to
the generators $e_i^*\in H^1(G,\bF_p)$,
because this set of generators %are 
is preserved, up to sign,
by the action of ${\rm Aut}(G)$.
In turn, it acts canonically on the spaces
of $1$-chains and $0$-chains
of ${\rm Alb}_p(M(G))$.

\begin{remark}
 The isomorphism \eqref{isom1} is not $\Aut(G)$-equivariant in general. The action of $a\in \Aut(G)$
 on %(non-oriented) 
 edges of ${\rm Alb}_p(M(G))$ sends
 $v \to v+e_i^*$ to the directed
 edge $a^*v \to a^*(v + e_i^*)$ 
 (whose direction may not agree with the one given 
 by the orientation of ${\rm Alb}_p(M(G))$),
 whereas the action on 
 $\bF_p[H^1(G,\bF_p)] \otimes C^1(G,\bF_p)$ 
 sends $[v]\otimes e_i^*$ to $ 
 [a^*v]\otimes (a^*(e_i^*))$. 
 Note that
$a^*( e_i^\ast) = \pm e_{a^{-1}(i)}^\ast$
with the sign depending on 
whether $a\colon a^{-1}(i)\mapsto i$ 
preserves or reverses orientation of the edge
$i\in E(G)$.
Hence the actions of $a\in {\rm Aut}(G)$
on the two sides of the isomorphism (\ref{isom1})
agree only when
$a$ preserves the orientation of $G$.
\end{remark}

By Proposition 
\ref{solution-rephrase}, 
solutions $b = (b_i)_{i\in E}$
of $G$ mod $p$ can be interpreted as elements 
of the morphism space
$\beta \in {\rm Hom}(C^1(G,\bF_p), 
\,C_1({\rm Alb}_p(M(G)),\bF_p))$. 
Then ${\rm Aut}(G)$ acts canonically on 
this morphism space via conjugation. 
That is, $a \cdot \beta = 
a \circ  \beta \circ a^{-1}$, for
$a \in \Aut(G)$.

\begin{proposition} \label{prop:solpres} 
Solutions of $M(G)$ mod $p$ 
form an ${\rm Aut}(G)$-invariant subspace
$${\rm Sol}_p(M(G))\subset 
{\rm Hom}(C^1(G,\bF_p), 
\,C_1({\rm Alb}_p(M(G)),\bF_p)).$$
Given $\beta\in {\rm Sol}_p(M(G))$,
the color profile
satisfies $\lambda \circ (a\cdot \beta)= 
a\cdot (\lambda \circ \beta)$ for any
$a\in {\rm Aut}(G)$.
In particular,
if $\beta$ has a constant color profile 
$\lambda\circ \beta = c\cdot {\rm Id}_{C^1(G,\,\bF_p)}$,
then $a\cdot \beta$ has a constant
color profile,
for the same constant $c\in \bF_p$.
\end{proposition}

\begin{proof}
Condition
(A) in Proposition \ref{solution-rephrase} is preserved,
because the action of $a\in {\rm Aut}(G)$ 
permutes the $E$-colorings of $G$ and 
${\rm Alb}_p(M(G))$ by the same permutation.
Thus, conjugation by $a$ respects colors,
i.e.~$(a\cdot \beta)(e_i^*) = (a\beta a^{-1})(e_i^*)$ 
is a $1$-chain
supported on the edges of ${\rm Alb}_p(M(G))$
of color $i$.
Condition (B) in Proposition \ref{solution-rephrase} 
is preserved, because the action 
of ${\rm Aut}(G)$ on the source and 
target of $\beta$ stabilizes the subspaces
$B^1(G,\bF_p)$ and  $H_1({\rm Alb}_p(M(G)),\bF_p)$, 
respectively. The first statement follows.

We have
$a\cdot (\lambda\circ \beta)  = 
(a\cdot \lambda)\circ (a\cdot \beta)$. Note that
$a\cdot \lambda = a\lambda a^{-1}=\lambda$ because
$\lambda$ is seen to be ${\rm Aut}(G)$-equivariant.
Thus, $\lambda\circ (a\cdot \beta) = a\cdot (\lambda\circ \beta)$.
Therefore, if 
$\lambda\circ \beta = c\cdot {\rm Id}_{C^1(G,\bF_p)}$, 
$$ a\cdot (\lambda\circ \beta) = 
a(c\cdot {\rm Id}_{C^1(G,\,\bF_p)}) 
a^{-1} = c\cdot {\rm Id}_{C^1(G,\,\bF_p)}.$$ 
So $\lambda \circ(a\cdot \beta) = 
\lambda\circ \beta$, i.e.\ $\beta$ and 
$a \cdot \beta=a\beta a^{-1}$
have the same constant color profile. 
\end{proof}

More generally,
the conjugation action permutes 
the color profile of a solution $\beta$, 
according to the underlying permutation of $a$.
We may unwind the above proposition in terms
of collections of $1$-cycles 
$ b= (b_i)_{i\in E}$ supported on the 
edges of color $i\in E$, as follows:

\begin{corollary}\label{unwind}
Let $b= (b_i)_{i\in E}$ be a 
solution of $M(G)$ mod $p$.
Then so is $b'\coloneqq (b_i')_{i\in E}$,
where $b_i'\coloneqq \pm a\cdot b_{a(i)}$
with sign determined by whether $a\in {\rm Aut}(G)$
respects or reverses the orientation
of the edge $i\mapsto a(i)$. If $b$ has constant
color profile, so does $b'$, for the same constant.
\end{corollary}

\begin{proof}
Suppose $b\leftrightarrow \beta$ under
the identification of 
Proposition \ref{solution-rephrase}. Let
$\beta' \coloneqq  a\cdot \beta = a\beta a^{-1}$.
On $C^1(G,\bF_p)$, 
the canonical action is
$a^{-1} \cdot e_i^* = \pm e_{a(i)}^*$ with
the sign determined as in the statement of the 
corollary.  By Proposition
\ref{prop:solpres}, 
$b_i'\coloneqq \beta'(e_i^*) = 
\pm a \beta(e_{a(i)}^*) = \pm a\cdot b_{a(i)}$
defines a solution $b'\leftrightarrow \beta'$. 
The last statement 
concerning color profiles also
follows from Proposition
\ref{prop:solpres}.
\end{proof}

The graph ${\rm Alb}_p(M(G))$ has 
an involution $\iota$, which
acts on vertices by the map
$$
\iota \colon H^1(G,\bF_p)\to 
H^1(G,\bF_p), \qquad v \mapsto -v,$$
and sends the edge $v\to v+ e_i^*$ to the edge
$-v \to -v-e_i^*$ (which is the edge
$-v-e_i^*\to -v$, but with the wrong orientation).
So $\iota$ respects
the colors of all edges, but also reverses the
orientations of all edges of ${\rm Alb}_p(M(G))$.
Define an action of $\iota$
on ${\rm Sol}_p(M(G))$ by 
\begin{align}
    \label{iota-def}
\iota\cdot \beta\coloneqq 
\iota_*\circ \beta \circ ({\rm negation}) = 
-\iota_*\circ \beta,
\end{align}
where negation acts on the source 
$C^1(G,\bF_p)$
of $\beta$.
It preserves solutions because the action of
$\iota$ respects colors
and because negation and the pushforward
map $\iota_*$ on $1$-chains both
preserve the homology $H_1(\Alb_p(M(G)), \bF_p)$.
From the action of $\iota$ on edges,
we see that $\lambda \circ (\iota\cdot \beta)=
\lambda \circ \beta$ and so the action (\ref{iota-def})
preserves the color profile. Furthermore, 
$\iota$ commutes
with the action of ${\rm Aut}(G)$ on the space
of solutions. 
To prove this, note that $\iota_*$ and 
negation commute with the actions
of $a\in {\rm Aut}(G)$ on the source 
and target of $\beta$, respectively.

\begin{definition}\label{symm} Let $A\subset 
{\rm Aut}(G)\times \bZ/2$ be any subgroup. 
The {\it symmetrization map}
is 
\begin{align}\label{average-v1}\begin{aligned} 
\sigma_A\colon 
{\rm Sol}_p(M(G))&\to {\rm Sol}_p(M(G)) \\
\beta& \mapsto \textstyle \sum_{(a,\,r)\,\in\, A} 
\iota^r\cdot a\cdot \beta. 
\end{aligned} 
\end{align} \end{definition}

\begin{corollary}\label{average-color}
If $\beta$ has a constant color profile,
then $\lambda\circ \sigma_A(\beta) = 
|A|\cdot (\lambda\circ \beta)$. 
Moreover, for any $(a,r) \in A$, 
we have $\iota^r\cdot a \cdot 
\sigma_A(\beta) = \sigma_A(\beta).$ 
\end{corollary}

In other words, the symmetrization map
sends a solution $\beta$ to one which
is $A$-invariant.
The corollary alleviates
some of the computational
difficulties involved in determining the existence
of a solution with constant nonzero color profile,
by reducing the solutions under consideration 
to a smaller space of invariant solutions. 

We now consider a special case,
relevant to the computations in the following section. 
Let $G$ be a graph, and choose an orientation. 
Let $A \subset \Aut(G) \times \bZ/2$ be a subgroup such that, 
for any element
$(a,r)\in A$, $a\in {\rm Aut}(G)$ preserves or reverses the 
orientation of every edge of $G$. Let 
\begin{align}\label{chi}\chi \colon 
A \to \set{\pm 1}\end{align} be the 
character satisfying
$\chi(e,1) = -1$ and $\chi(a,0) = 1$ if and only 
if $a \in \Aut(G)$ preserves the orientation of 
every edge of $G$. 
\begin{definition}
With the above notation and assumptions, we say that 
a $1$-chain $b \in C_1(\Alb_p(M(G)), \bF_p)$ is 
\emph{$(A,\chi)$-equivariant} if 
$\iota^r \cdot a \cdot b = \chi(a,r)b$ 
for every $(a,r) \in A$. 
\end{definition}

\begin{proposition}\label{prop}
Let $G$ be an oriented graph and 
$A \subset \Aut(G) \times \bZ/2$ be a subgroup such 
that for any $(a,r) \in A$, $a$ either preserves 
or reverses the orientation of every edge of $G$.
Then, on the level of $1$-chains 
$b = \sum_{e\in E({\rm Alb}_p(M(G)))} b_e e$ 
of $\Alb_p(M(G))$, the 
symmetrization map is 
\begin{align}\label{average-v1-1chains}\begin{aligned}
\sigma_A\colon 
C_1(\Alb_p(M(G)),\bF_p)&\to C_1(\Alb_p(M(G)),\bF_p),\\
b & \mapsto \textstyle \sum_{(a,\,r)\,\in\, A} 
\chi(a,r)
(\iota^r\cdot a\cdot b). 
\end{aligned} 
\end{align}
The element $\sigma_A(b)$ is a solution whenever $b \in C_1(\Alb_p(M(G)), \bF_p)$ is a solution. 
Its color profile is $\lambda(\sigma_A(b)) = 
\va{A} \cdot \lambda(b)$ 
if the color profile of $b$ is constant.

Finally, when $|A|$ and $p$ are coprime,
the image of the map $\sigma_A$ in 
\eqref{average-v1-1chains} is the space of 
$(A,\chi)$-equivariant $1$-chains
$b \in C_1(\Alb_p(M(G)),\bF_p)$. 
\end{proposition}

\begin{proof}
The fact that $\sigma_A$ preserves solutions 
follows from Corollary \ref{unwind} 
and the
ensuing discussion of $\iota$. The discrepancy
in signs between (\ref{average-v1}) and 
(\ref{average-v1-1chains}) is due to the additional
global negation arising from Corollary \ref{unwind},
when an automorphism reverses all edge orientations,
and the sign appearing in the definition 
(\ref{iota-def}) 
of the action of $\iota$.
Corollary \ref{average-color} implies
the color profile multiplies by $|A|$ when
the color profile of $b$ is constant. For 
the last statement, note that the averaging 
$\sigma_A$
produces a $1$-chain
which is $A$-equivariant with respect 
to the character $\chi$
and that conversely, $\sigma_A(|A|^{-1}b) = b$
for any $(A, \chi)$-equivariant $1$-chain 
$b\in C_1(\Alb_p(M(G)), \bF_p)$.
\end{proof}

\begin{remark}
The
symmetrization map 
can alternatively be described as follows.  
For a positively oriented edge 
$e \in E(G)$, put
\begin{align} \label{absval}
\va{\pm e} \coloneqq 
e \in C_1(G,\bF_p) = \bF_p^{E(G)}. 
\end{align}
That is, $\va{e} = e$ and $\va{-e} = e \in C_1(G,\bF_p)$.
This ``absolute value'' is only 
defined for $1$-chains of the form $\pm e\in C_1(G,\bF_p)$. 
Let $A \subset \Aut(G) \times \bZ/2$ be a subgroup 
such that, for $(a,r) \in A$, $a$ preserves or reverses 
the orientation of every edge of $G$. Then the 
symmetrization map sends an edge 
$e \in C_1(\Alb_p(M(G)), \bF_p)$ 
to the element $\textstyle 
\sum_{(a,\,r)\,\in\, A} 
|\iota^r\cdot a\cdot e|$.
\end{remark}

\section{Proof of the main result} \label{sec:proof}

Here we prove the following theorem.

\begin{theorem}\label{average-no} 
Let $K_7$ be the complete graph
on $7$ vertices. Let $A = D_7\times 
\bZ/2$
where $D_7\subset {\rm Aut}(K_7)$ is the standard
dihedral action. Then every $A$-invariant 
element 
$\beta\in {\rm Sol}_3(M(K_7))$
with constant color profile
satisfies $\lambda\circ \beta=0$.
\end{theorem}

The proof will be reduced to 
a computer calculation, detailed in 
Section \ref{sec:comp} below.
The code is available
in the repository \cite{code}. 
We now give the theoretical
framework underpinning the computer computation.

\begin{figure}
\includegraphics[width = 2.9in]{tournament.png}
\caption{Orientation of $K_7$ as a 
regular tournament with
step-1 edges (blue), step-2 edges (green), and
step-3 edges (yellow).}
\label{tournament}
\end{figure}

\subsection{Lexicographic vertex fundamental domain}
\label{sec-fund}
We set up the computation in the generality of a
complete graph $K_{2n+1}$ and an arbitrary
prime $p\nmid 2n+1$. It is convenient to label
the vertices by $\{0,1,\dots, 2n-1, 2n\} = \bZ/(2n+1)$ 
and the oriented edges as 
$$E(K_{2n+1}) = \{e_{i,\,i+1}, \,e_{i,\,i+2}, 
\, e_{i,\,i+3}, \,\dots,\, 
e_{i, \,i+n}\}_{i\textrm{ mod }2n+1}.$$ 
That is, we orient
$K_{2n+1}$ as a regular tournament, see
Figure \ref{tournament} for the example
$2n+1=7$. We call $e_{i,\,i+k}$ the 
{\it step-$k$ edges}. We require:

\begin{lemma}\label{step-star}
For any prime $p\nmid 2n+1$, 
the linear functionals 
$\{e_{i,\,i+k}^*\}_{i = 0,\, \dotsc, \,2n, \,k=1,\,\dots\,,\,n-1}\cup \{s^*\}$
form a basis of $H^1(K_{2n+1},\bF_p)$, where $s^* \coloneqq 
\sum_{i\,\,{\rm mod }\,\,2n+1} e_{i,\,i+n}^*$ is the 
{\rm star dual}.
\end{lemma}

\begin{proof}
Observe that $h^1(K_{2n+1}, \bF_p) = 
\va{E(K_{2n+1})}-\va{V(K_{2n+1})} + 1 = 
{2n+1\choose 2} - (2n+1) + 1$ 
which agrees with the number of linear
functionals provided. Thus, it suffices to check
that any cycle $\gamma\in H_1(K_{2n+1},\bF_p)$ on which
all the linear functional vanish, is necessarily $\gamma=0$.

The vanishing of the step-$k$ edge duals, 
$e_{i,\,i+k}^*(\gamma)=0$, $k\in \{1,\dots,n-1\}$,
implies that $\gamma$
is a cycle supported on the step-$n$ edges. Up
to taking a multiple, 
there is a unique such cycle, which
is the sharpest $(2n+1)$-pointed star
$s
= \sum_{i\,\,{\rm mod }\,\,2n+1} e_{i,\,i+n}
\in H_1(K_{2n+1},\bF_p)$, depicted
in Figure \ref{tournament} as the sum
of the yellow edges. Hence $\gamma = m \cdot s$ 
for some $m \in \bF_p$. 
But then we have 
$s^\ast(s)$
$= 2n+1 \not\equiv 0\textrm{ mod }p$.
Thus the additional condition
$s^*(\gamma)=0$ ensures $\gamma=0$.
\end{proof}

We henceforth let $n \in \bN$, and let $p$ be 
a prime number that does not divide $2n+1$.
We call the linear functionals of Lemma \ref{step-star}
the {\it step-star basis} of $H^1(K_{2n+1},\bF_p)$.
The major advantage of the step-star basis
is that it is naturally 
$D_{2n+1}$-equivariant. More precisely, given an element 
$a\in C_{2n+1}\subset D_{2n+1}$ of the cyclic subgroup, we have
$a\cdot e_{i,\,i+k}^* = e_{i-a,\,i+k-a}^*$,
and given an element $a\in D_{2n+1}\setminus C_{2n+1}$,
the action of $a$ also preserves the set of step-$k$ edge
duals and the star dual, after a global negation 
(such an $a$ visibly reverses the orientation of every
edge of the chosen orientation of $K_{2n+1}$,
see Figure \ref{tournament}).

To set up the computation explicitly, we must
find a fundamental domain for the action of 
$D_{2n+1}\times \bZ/2$ (or more generally, 
a subgroup $A\subset D_{2n+1}\times \bZ/2$), on 
$V({\rm Alb}_p(M(K_{2n+1}))) = H^1(K_{2n+1},\bF_p)$
and elements of 
$E({\rm Alb}_p(M(K_{2n+1})))$
representing every $A$-orbit.

The {\it lexicographic fundamental domain}  
$V_{\rm fund}\subset V({\rm Alb}_p(M(K_{2n+1}))) 
= H^1(K_{2n+1},\bF_p)$
is defined as follows. For $v\in H^1(K_{2n+1},\bF_p)$, 
take the unique
expression of $v$
in the step-star basis $$ 
v = cs^*+ \sum_{k=1}^{n-1} \sum_{i=0}^{2n}
c_{i,\,i+k} e_{i,\,i+k}^* \qquad \textrm{ with } 
\qquad c, c_{i,\,i+k}\in \{0,1,\dots,p-1\}.$$
Consider now the action of 
$(a,r)\in D_{2n+1}\times \bZ/2$
on $v$. We begin by taking the lexicographically minimal
$D_{2n+1}\times \bZ/2$ orbit-representative
of the step-$1$ coefficients
$(c_{01}, c_{12}, c_{23}, \dots, c_{2n,0})$
with respect to the ordering $0<1<\cdots< p-1$.
This amounts to choosing the minimal representative
of such a $\{0, \dots, p-1\}$-string under
the operations of cyclic rotation, reflections,
and negations. Having found this minimal 
representative $c_{\textrm{step-}1}$,
we set $A':={\rm Stab}_{A}(c_{\textrm{step-}1})$.

Then we consider the $A'$-action on the step-$2$ 
coefficients
$(c_{02}, c_{13}, c_{24}, \dots, c_{2n-1,0}, c_{2n,1})$.
We again find a lexicographically minimal 
representative under the $A'$-action 
$c_{\textrm{step-}2}$, and define $A'':={\rm Stab}_{A'}(c_{\textrm{step-}2})$. Continuing inductively, 
we arrive at a subgroup $A^{(n)}\subset A$. It acts
on the coefficient $c$ of the star dual element
$s^*$ by the character
$(a,r)\cdot c = (-1)^{{\rm sgn}(a)+r}c.$ 
Here, $\rm{sgn}$ is the orientation 
character of the dihedral group,
which preserves or reverses the
orientation of every edge of $K_{2n+1}$.
Again taking the lexicographically minimal
representative for $c$, we reduce any vertex 
$v\in V({\rm Alb}_p(M(K_{2n+1})))$ 
into a fundamental domain, expressible in the form 
$$v_{\rm fund}:=(c_{\textrm{step-}1}, 
c_{\textrm{step-}2},\dots, c_{\textrm{step-}(n-1)}, 
c_{\rm star}).$$
We additionally extract the stabilizer group
${\rm Stab}_A(v_{\rm fund}) = {\rm Stab}_{A^{(n)}}(c_{\rm star})$.

\begin{definition}
    Define $V_{\rm fund} \subset  H^1(K_{2n+1},\bF_p)$, 
    the \emph{lexicographic
    fundamental domain} of the vertices
    under the action of $A$, %consisting
    as the set of all $v_{\rm{fund}}$ for $v \in H^1(K_{2n+1},\bF_p)$. 
\end{definition}

\begin{example}
    For ${\rm Alb}_3(M(K_7))$, 
    consider the vertex represented by 
    $v = (1200220, 2200111, 1)$, with ternary
    strings representing the step-$1$, step-$2$,
    and star coefficients, respectively. Then 
    we have $v_{\rm fund} = (0011021, 0022211, 2)$
    with a group element 
    satisfying $a\cdot v = v_{\rm fund}$ given by 
    $a = ({\rm rot}^2, 1)\in D_7\times \bZ/2 $ and 
    ${\rm rot}\in C_7$ the cyclic rotation increasing
    the indices of vertices by $1$.
    The stabilizer ${\rm Stab}_{D_7\times \bZ/2}
    (v_{\rm fund})$ is trivial.
\end{example}

\subsection{Equivariant $1$-chains of 
the Albanese graph}
Let $v=v_{\rm fund}\in V_{\rm fund}$ be a fundamental
domain vertex. We consider, for each $i\textrm{ mod }2n+1$,
the incoming and outgoing edges of $v$ 
of color $i$, given respectively by 
\begin{align*}
e_{\rm in}(v, i) &
\coloneqq 
(v-e_i^*\to v) \quad
\textrm{ and } \quad
e_{\rm out}(v, i)\coloneqq 
(v\to v+ e_i^*)\in 
E({\rm Alb}_p(M(K_{2n+1}))).\end{align*}
For each ingoing or outgoing edge $e$ from 
$v\in V_{\rm fund}$, 
we take the $(A,\chi)$-average 
\begin{align}\label{average-v2}
\sigma_A(e)\coloneqq
\sum_{(a,r)\in A} |\iota^r \cdot a \cdot e|
\in C_1({\rm Alb}_p(M(K_{2n+1})), \bF_p),
\end{align}
where the ``absolute value'' $|\cdot|$ 
is defined by 
\eqref{absval}.
We define the 
{\it color $j$ boundary} of $\sigma_A(e)$ as follows. 
Break $\sigma_A(e)$
into its color components
$$\sigma_A(e) = \sum_{j\in E(K_{2n+1})}
\sigma_A(e)_j$$
with $\sigma_A(e)_j\in C_1({\rm Alb}_p(M(K_{2n+1})),\bF_p)$ 
supported on edges of color $j$, 
and take the boundary 
$$\partial \sigma_A(e)_j \in 
C_0({\rm Alb}_p(M(K_{2n+1}))) = 
\bF_p^{V({\rm Alb}_p(M(K_{2n+1})))}.$$
Define a {\it color weighting} as an element
$w\in \bF_p^{E(K_{2n+1})} = C^1(K_{2n+1},\bF_p)$. 
For a given
color weighting $w$, define the {\it $w$-weighted
boundary} by the following formula: 
$$\partial_{A,w}(e)
\coloneqq \sum_{j\in E(K_{2n+1})} w_j 
\partial_A(e)_j \vert_{V_{\rm fund}}\in 
\bF_p^{V_{\rm fund}}.$$
That is, $\partial_{A,w}(e)$ represents the boundary
of an
averaged $1$-chain (with respect
to the action (\ref{average-v2})), which
has been further weighted by $w_i$ on the $i$-th color,
and then restricted to 
 $V_{\rm fund}$. 
 
 Consider the character 
$\chi \colon A \to \set{\pm 1}$, see equation \eqref{chi}.

\begin{lemma}\label{sol-ok}
Let $w^{(1)}, \dotsc, w^{(2n)}$ be a basis of $
B^1(K_{2n+1},\bF_p)$.
A $1$-chain
 \begin{align}\label{form}b = \sum_{v\in V_{\rm fund}}
 \,\,\sum_{i\,\in\, E(K_{2n+1})}\,\,
 \sum_{\substack{e\,=\,e_{\rm in}(v, \,i)\\ 
 \textrm{or }\,e_{\rm out}(v, \,i)}} 
 c_e \sigma_A(e)\in C_1({\rm Alb}_p(M(K_{2n+1})), \bF_p)\end{align}
with $c_e\in \bF_p$
is a solution if and only if
$\partial_{A,w^{(k)}}(b)=0$ for all $w^{(k)}$.   
When $|A|$ and $p$ are coprime,
every $(A, \chi)$-equivariant
$1$-chain on ${\rm Alb}_p(M(G))$
is of the form (\ref{form}).
\end{lemma}

\begin{proof}
The condition that $\partial_{A,w^{(k)}}(b)=0$ is equivalent
to the condition that $\sum_{i\in E(K_{2n+1})}
w^{(k)}_i b_i\in H_1({\rm Alb}_p(M(G)),\bF_p)$ is a 
$1$-cycle, because for an 
equivariant $1$-chain,
it suffices to check
that the boundary vanishes on an $A$-fundamental
domain of vertices. 
Thus, if this condition holds for $w^{(k)}$ 
ranging over a basis
of $B^1(K_{2n+1},\bF_p)$, then $b$ defines a solution.

If $|A|$ and $p$ are coprime, then
any $(A, \chi)$-equivariant $1$-chain $b$
is expressible as a linear combination
of $\sigma_A(e)$, for $e$ ranging over 
$e_{\rm in}(v,i)$, $e_{\rm out}(v,i)$ for all
$v\in V_{\rm fund}$ and $i\in E$. 
This holds by Proposition \ref{prop},
as every edge of ${\rm Alb}_p(M(G))$
can be reduced (possibly
non-uniquely) via $A$, to an edge incident to a 
fundamental domain vertex.
Hence, $b$ is of the form \eqref{form}. 
\end{proof}

\subsection{Building the equivariant solution matrix}
\label{sec-sol}

We now define a matrix whose kernel admits a natural 
surjection onto the space of $(A, \chi)$-equivariant 
solutions $b \in C_1(\Alb_p(M(G)), \bF_p)$.

\begin{definition}
Define a matrix $S_A(p, 2n+1)$ over $\bF_p$, as follows. 
It is a vertical stacking of ${\rm rank}\,B^1(K_{2n+1},\bF_p)=2n$ matrices, each one
having proportions $ 
|V_{\rm fund}|\times (2n+1)(2n)\cdot |V_{\rm fund}|$. 
The column of the $k$-th matrix 
in the vertical stack corresponding to an ingoing
or outgoing edge $e$ to $v\in V_{\rm fund}$ is given
by the entries of 
$\partial_{A,w^{(k)}}(e)\in \bF_p^{V_{\rm fund}}$. 
The result is a matrix with
proportions $(2n)\cdot |V_{\rm fund}|\times 
(4n^2+2n)\cdot |V_{\rm fund}|$,
that we call $S_{A}(p,2n+1)$. 
\end{definition}

\begin{remark}
It follows from Lemma \ref{sol-ok} that the right
kernel of $S_A(p,2n+1)$ is the space of expressions
of the form (\ref{form}) which represent 
an $A$-invariant 
(in $\beta$ form) element
of ${\rm Sol}_p(M(K_{2n+1}))$. In particular,
$\ker(S_A(p,2n+1))$ 
surjects onto
${\rm Sol}_p(M(K_{2n+1}))^A$. But, in general, 
there may be a kernel, because the columns of 
$S_A(p,2n+1)$ do not represent a 
basis of the space of 
$(A,\chi)$-equivariant
$1$-chains (in $b$ form) on 
${\rm Alb}_p(M(K_{2n+1}))$. The issue is 
that
$e_{\rm in}(v, i)$, 
$e_{\rm out}(v, i)$ 
for $v\in V_{\rm fund}$ and $i\in E(K_{2n+1})$
do not form a fundamental domain for the $A$-action
on $E({\rm Alb}_p(M(K_{2n+1})))$, 
though they do contain all $A$-orbit representatives.
The averaging $\sigma_A(e)$ will identify
averaged $1$-chains associated to $e$, $e'$
if and only if they lie in the same $A$-orbit. 
As such, we only know 
$\dim \,\ker(S_A(p,2n+1))\geq 
\dim\,{\rm Sol}_p(M(K_{2n+1}))$. 
\end{remark}

The condition that an $(A, \chi)$-invariant
solution $b$ has color profile zero can be
detected from an expression of
the form (\ref{form}).
To detect such solutions,
define a matrix $S_A^{\textrm{step-}1}(p, 2n+1)$
by appending a single row to 
$S_A(p, 2n+1)$ which has an entry $1$ on all
columns indexed by 
$e_{\rm in}(v, i)$, 
$e_{\rm out}(v, i)$ 
for $i\in E(K_{2n+1})$ a step-1 edge,
and a $0$ entry for all other columns.
Similarly, we define 
$S_A^{\textrm{step-}1,2}(p, 2n+1)$
by appending to $S_A^{\textrm{step-}1}(p, 2n+1)$
a single row which is an indicator
on all step-$2$ edges, etc.

\begin{lemma}\label{stepjump}
If all $(A,\chi)$-equivariant solutions
are colorless, we 
have $${\rm rank}\,
S_A(p, 2n+1) ={\rm rank}\,S_A^{{\rm step}\textrm{-}1}(p, 2n+1) 
= S_A^{{\rm step}\textrm{-}1,2}(p, 2n+1)=\cdots $$
and conversely, 
when $|A|$ and $p$ are coprime,
and there is an $(A, \chi)$-equivariant solution 
with constant nonzero color profile, then 
$ {\rm rank}\,S_A^{\textnormal{step-}1}(p, 2n+1)= 
{\rm rank}\,
S_A(p, 2n+1) +1.$
\end{lemma}
\begin{proof}
Indeed,  $\ker(S_A^{\textrm{step-}1}(p, 2n+1))$ surjects
onto the space of
$(A,\chi)$-equivariant solutions 
which have color profile zero on all step-$1$ edges,
$\ker(S_A^{\textrm{step-}1,2}(p, 2n+1))$ surjects
onto the space of 
$(A,\chi)$-equivariant solutions 
which have color profile $0$ on 
step-$1$ and $2$ edges,
etc. (Here we used that the number of step-$k$ edges
is $2n+1$, which is coprime to $p$.)

The converse follows from Proposition \ref{prop}.
\end{proof}

This concludes the theoretical setup for the 
computation proving Theorem \ref{average-no}.

\subsection{Computational aspects of the proof of
Theorem \ref{average-no}} \label{sec:comp}
Following the logic of the computation
laid out in Sections \ref{sec-fund}--\ref{sec-sol},
the code
{\tt SymmetrizedSolutionsKn}
\cite{code} executes the following steps: 
\begin{enumerate}
\item[($V_{\rm fund}$)] Build the vertices 
$ V({\rm Alb}_p(K_{2n+1}))$ 
using the step-star basis, and extract
the lexicographic fundamental domain 
$V_{\rm fund}\subset V({\rm Alb}_p(K_{2n+1}))$
with respect to the action of 
$A\subset D_{2n+1}\times \bZ/2$ labeled as 
$A ={\tt DnZ2,\, Dn,\, CnZ2,\, Cn,\, Z2, \,Triv}$.
Nota bene: The program uses $n$ for $2n+1$. \smallskip
\item[($E_{\rm orb}$)] Build the ingoing
and outgoing directed edges 
$v-e_i^*\to v$ and
$v\to v+e_i^*$ by ranging over all
$v\in V_{\rm fund}$ in a fundamental
domain and all colors $i$.
\smallskip
\item[$(w^{(k)})$] Choose a basis 
$(w^{(k)})_{k=1,\,\dots,\,2n}$ of the space 
$B^1(K_{2n+1},\bF_p)$.
\smallskip
\item[($S_A(w^{(k)})$)] For each color
weighting $w^{(k)}$, 
build the $w^{(k)}$-weighted coboundary matrix,
which is the map 
$S_A(w^{(k)})\colon \bF_p^{E_{\rm orb}} 
\to \bF_p^{V_{\rm fund}}$ whose kernel surjects
onto the $(A,\chi)$-equivariant
$1$-chains $b\in C_1({\rm Alb}_p(M(K_{2n+1})),\bF_p)$ 
for which  $\sum_{i\in E(K_{2n+1})} w_i^{(k)}b_i$ 
are $1$-cycles. 
\smallskip
\item[($S_A$)] Stack the matrices $S_A(w^{(k)})$, ranging
over the basis $w^{(k)}$, $k=1,\,\dots,\,2n$.
\smallskip
\item[($S_A^{\textrm{step-}1,\dots,s}$)] 
For $s=1,\dots,n$
add, for each step size $i\in \{1,\dots,s\}$, 
a row to $S_A$, which is an indicator function 
on the edges of step size $i$.
\smallskip
\item[($r$, $r^s$)] Compute the ranks 
$r\coloneqq {\rm rank}(S_A)$, $r^s 
\coloneqq {\rm rank}(S_A^{\textrm{step-}1,\dots,s})$, 
using sparse rank functions; run by
{\tt TestColorlessAveragedSolutions(p, 2n+1~:~subgroup:="}$A${\tt ")}.
\smallskip
\item[(answer)] If $r=r^1$, 
all $(A, \chi)$-equivariant solutions with constant color profile
have zero color profile, see Lemma \ref{stepjump}.
\end{enumerate}

Strictly speaking, we only need the computation
up to the step $s=1$; the higher values are built-in as
a check. The computation is done in MAGMA \cite{magma}
which is well-suited to rank computations of 
very sparse, large matrices over finite fields. 

\begin{proof}[Proof of Theorem \ref{average-no}] 
Performing the above computation for $p=3$ and $2n+1=7$,
and $A= D_7\times \bZ/2$ (which indeed satisfies
that $|A|=28$ is coprime to $3$ and that $p$ and $2n+1$
are coprime), the outcomes are as follows:
The size of the fundamental
domain is $|V_{\rm fund}| = 517{,}570$. 
There are 
${7 \choose 2}=21$ colors, and the number of color
weighting vectors $w^{(k)}$ is 
${\rm rank}\,B^1(K_7,\bF_3) = |V(K_7)|-1=6$.
We get that $S_A(3,7)$ is a 
$3{,}105{,}420\times 21{,}737{,}940$ matrix over 
$\bF_3$. As expected, this is
$(2n)\cdot |V_{\rm fund}|\times 
(2n+1)(2n)\cdot |V_{\rm fund}|$.
The matrices, $S_A^{\textrm{step-}1}(3,7)$,
$S_A^{\textrm{step-}1,2}(3,7)$,
$S_A^{\textrm{step-}1,2,3}(3,7)$ respectively
have $1$, $2$, $3$ more rows. 
We find that the ranks 
$r = r^{1} = r^{2} = r^{3}= 3{,}060{,}117$
are all equal, though we only need
the first equality for a proof. 
The total runtime of the program
was $132{,}396$ seconds ($\approx 37$ hours),
executed on a Gemini HPC node of Utrecht
University with 32 CPU cores and 8 GB of 
memory per core. The peak memory usage was 114 GB.
\end{proof}

\begin{corollary}\label{cor} Any element of
${\rm Sol}_3(M(K_7))$ with constant color
profile has zero color profile.\end{corollary}

\begin{proof} By Corollary \ref{average-color},
any solution with constant nonzero
color profile can be
$A=D_7\times \bZ/2$-averaged
to a solution with constant 
nonzero color profile, since $|A|$ and $3$
are coprime. We deduce
the corollary from Theorem \ref{average-no}. \end{proof}

\begin{remark}\label{computation-size}
One may wonder; what is the point in the $A$-averaging 
and reduction into a fundamental domain? The reason
is that the rank computations are the rate-limiting
step of the computation. Without symmetrizing at all,
the relevant matrix has proportions 
$$(2n)\cdot p^{2n^2-n} \times 
(2n+1)(n)\cdot p^{2n^2-n},$$
where $p^{2n^2-n} = 
|V({\rm Alb}_p(M(K_{2n+1})))|$. 
(Whenever $A$ preserves the orientations of all 
edges of ${\rm Alb}_p(M(K_{2n+1}))$, the
outgoing edges alone from each 
fundamental domain vertex represent
every $A$-orbit of edges. So the number of columns
can be halved.) It 
has a similar sparsity to $S_A(p, 2n+1)$.
For $p=3$ and $2n+1=7$, we 
get a matrix over $\bF_3$ with proportions
$$86{,}093{,}442 \times 301{,}327{,}047.$$
Finding the rank of such a matrix was beyond
our computing capacity.
\end{remark}

\begin{proof}[Proof of Theorem \ref{main}]
It follows from \cite[Cor.~2]{ADFIO} that, 
for any $(X,\Theta)\in \cA_6$, the class
$6\cdot \Theta^5/5!$ is represented by an  
effective curve, arising from
a Prym--Tyurin construction of exponent $6$.
The original construction is due to Kanev 
\cite{kanev} but dominance onto $\mathcal A_6$ was proven in \cite{ADFIO}. 
Thus, the minimal multiple of the minimal
class which is represented
by an algebraic cycle always divides $6$.
Let $K_7$ be the complete graph on $7$ vertices.
By the results of \cite{companion},
see Theorem \ref{thm:divides}, 
the combinatorial
invariant $d(K_7)\in \bN$ divides $N_6$. 
It follows from Corollary \ref{cor}
that $3\,{\mid}\, d(K_7)$. By Proposition
\ref{div}, we have $d(K_5)
\,{\mid}\, d(K_7)$. Finally, $2\,{\mid}\, d(K_5)$ by 
\cite[Prop.~7.5]{companion}. %So $d(K_7)=6$ 
So $6\mid d(K_7) \mid N_6 \mid 6$
and hence $N_6=6$.
\end{proof}

\printbibliography

\end{document}